\documentclass[10pt, article, titlepage,openany, rotate]{amsart}

\usepackage{ae} 
\usepackage[T1]{fontenc}
\usepackage[cp1250]{inputenc}
\usepackage{amsmath}
\usepackage{amssymb, amsfonts,amscd,verbatim}
\usepackage{color}
\usepackage[colorlinks]{hyperref}
\usepackage[dvips]{graphicx}
\usepackage{latexsym}
\usepackage{eucal,eufrak}
\usepackage{indentfirst}
\usepackage{latexsym}

\usepackage{graphicx}
\usepackage{verbatim}   
\usepackage{color}      
\usepackage{subfigure}  

\usepackage{amsmath}    

\theoremstyle{plain}
\newtheorem{Prop}{Proposition}[section]
\newtheorem{Thm}[Prop]{Theorem}
\newtheorem{Cor}[Prop]{Corollary}

\newtheorem{Lem}[Prop]{Lemma}

\theoremstyle{definition}
\newtheorem{Def}[Prop]{Definition}

\theoremstyle{remark}
\newtheorem{Rem}[Prop]{Remark}
\newtheorem{Problem}[Prop]{\bf Problem}

\def\dim{\mathop{\roman{dim}}}
\def\int{\mathop{\roman{int}}}

\def\1{^{-1}}

\def\N{{\bold N}}
\def\Z{{\bold Z}}
\def\Q{{\bold Q}}

\def\PP{{\mathbb P}}

\def\dim{\text{dim}}

\def\Tor{\text{Tor}}

\def\Ab{\text{Ab}}

\def\PPP{{\mathcal P}}
\def\NNN{{\mathcal N}}

\def\dokaz{{\bf Proof. }}
\numberwithin{equation}{section}

\input pstricks.tex
\input xy
\xyoption{all}

\begin{document}
\title[
BOCKSTEIN THEOREM FOR NILPOTENT GROUPS
]%
   {BOCKSTEIN THEOREM FOR NILPOTENT GROUPS}

\author{M.~Cencelj}
\address{IMFM,
Jadranska ulica 19,
SI-1111 Ljubljana,
Slovenija }
\email{matija.cencelj@guest.arnes.si}

\author{J.~Dydak}
\address{University of Tennessee, Knoxville, TN 37996, USA}
\email{dydak@math.utk.edu}

\author{A.~Mitra}
\address{University of Tennessee, Knoxville, TN 37996, USA}
\email{mitra@math.utk.edu}

\author{A.~Vavpeti\v c}
\address{Fakulteta za Matematiko in Fiziko,
Univerza v Ljubljani,
Jadranska ulica 19,
SI-1111 Ljubljana,
Slovenija }

\email{ales.vavpetic@fmf.uni-lj.si}

\date{May 17, 2008}


\thanks{Supported in part by the Slovenian-USA research grant BI--US/05-06/002 and the ARRS
research project No. J1--6128--0101--04}
\thanks{The second-named author was partially supported byĘĘMEC, MTM2006-0825.}

\begin{abstract}

We extend the definition of Bockstein basis $\sigma(G)$ to nilpotent
groups $G$. A metrizable space $X$ is called a {\it Bockstein space} if $\dim_G(X) =
\sup\{\dim_H(X) | H\in \sigma(G)\}$ for all Abelian groups $G$. Bockstein First Theorem
says that all compact spaces are Bockstein spaces.

Here are the main results of the paper:
\begin{Thm}
Let $X$ be a Bockstein space.
If $G$ is nilpotent,
then $\dim_G(X) \leq 1$ if and only if $\sup\{\dim_H(X) | H\in\sigma(G)\}\leq 1$.
\end{Thm}

\begin{Thm}
$X$ is a Bockstein space if and only if $\dim_{\Z_{(l)}} (X) = \dim_{\hat{Z}_{(l)}}(X)$
for all subsets $l$ of prime numbers.
\end{Thm}

\end{abstract}

\maketitle

\medskip
\medskip
\tableofcontents

\section{Introduction}

We use the Kuratowski notation $X\tau M$ in the case when every map from a closed 
subset of $X$ to $M$ can be extended over all $X$.

Recall that the cohomological dimension of a space $X$ with respect to Abelian group $G$
is less than or equal $n$, denoted by $\dim_G X\le n$, if $\check{H}^{n+1}(X,A;G)=0$
for all closed $A\subset X$. 

Of basic importance in cohomological dimension theory is the {\it Bockstein basis}
$\sigma(G)$ of Abelian group $G$ (see \ref{DefOfMaxBocksteinBasis})
and the following result of Bockstein (see \cite{K} or \cite{DR4}):
\begin{Thm}[Bockstein First Theorem] \label{BocksteinFirstThm}
If $X$ is a compact space, then $$\dim_G(X)=\sup\{\dim_H(X) | H\in\sigma(G)\}.$$
\end{Thm}

The aim of this paper is to generalize \ref{BocksteinFirstThm} to nilpotent groups.
There are two issues to resolve first:
\begin{enumerate}
\item Define cohomological dimension with respect to non-Abelian groups.
\item Define the Bockstein basis of nilpotent groups.
\end{enumerate}

The definition of $\dim_G(X)$ for $G$ non-Abelian was first introduced by
A.Dranishnikov and D.Repov\v s \cite{DrR1} as follows:
By \cite[Theorem 1.1]{DR4}, $\dim_G(X)\leq n$ is equivalent (for Abelian groups $G$) to 
$X\tau K(G,n)$, where $K(G,n)$ is an Eilenberg-MacLane space.
One ought to use the same equivalence in definition of $\dim_G(X)$ for non-Abelian groups.
For nilpotent groups $G$ this definition was used for the characterisation of nilpotent
CW complexes as absolut extensiors of metric compacta (see \cite{CD}).
There is no
Eilenberg-MacLane space $K(G,n)$, $1 < n <\infty$, for non-Abelian groups $G$, so 
$\dim_G X\in\{0,1,\infty\}$. Since $\dim_G X=0$ means $\dim(X)=0$ (see Lemma~\ref{0dim})
the only interesting question is if $X\tau K(G,1)$ holds.
\par Our definition of the Bockstein basis for nilpotent groups can be found in \ref{BocksteinBasisForNil}.
\par The remainder of this section is devoted to elementary properties
of cohomological dimension over non-Abelian groups.

\begin{Lem} \label{0dim}
If $X$ is a metrizable space and $\dim_G(X)=0$ for some
group $G\ne 1$, then $\dim_H(X)=0$ for any group $H$.
\end{Lem}
\dokaz
We will show that for nontrivial group $G$ we have $\dim_G(X)=0$ iff
$\dim_{\Z_2}(X)=0$  and $\dim_{\Z_2}(X)=0$ means $\dim(X)=0$..

Suppose $\dim_G(X)=0$ for some $G\ne 1$. Let $i\colon \Z_2\to G$ be an injection
and let $r\colon G\to \Z_2$ be a map such that $r\circ i=id_{\Z_2}$. Let $A$
be a closed subset
of $X$ and $f\colon A\to K(\Z_2,0)=\Z_2$ a map. Then there exists an
extension
$f'\colon X\to G$ of $i\circ f\colon A\to G$. Then $F=r\circ f'\colon X\to
\Z_2$
is an extension of $f$.

Suppose that $\dim_{\Z_2}(X)=0$.
Let $A\subset X$ be a closed subset and let $f\colon A\to G$ be continuous
map.
For every $g\in G$ we define $A_g=f^{-1}(g)$ and $B_g=\{x\in X\mid
d(x,A\setminus A_g)\le d(x,A_g)\}$.
Because $A\setminus A_g$ is a closed subset of $X$, $B_g$ is a closed subset of
$X$. Because $B_g\cap A_g=\emptyset$
we can define a continuous map $f_g\colon A_g\cup B_g\to \Z_2$ by
$f_g(A_g)=1$ and $f_g(B_g)=0$.
Let $F'_g\colon X\to \Z_2$ be an extension of $f_g$ and $F_g\colon X\to G$
be the composition of $F'_g$ and the inclusion $i_g\colon \Z_2\to G$ which
maps $i_g(0)=e$ and $i_g(1)=g$.
For every $g\in G$ we define $X_g=F_g^{-1}(g)$. Sets $X_g$ are open and
close in $X$ and
pairwise disjoint. So we can define a continuous map $F\colon X\to G$ as
$$
F(x)=\left\{
\begin{array}{ll}
g&; x\in X_g\\
e &; x\not\in \cup_{g\in G} X_g
\end{array}\right.
$$
The map $F$ is an extension of $f$.
\hfill
$\blacksquare$

\begin{Lem} \label{BocksteinInequalities}
Let $X$ be a metrizable space.
If $1\to J\to G\to I\to 1$ is an exact sequence of groups
and $\dim_JX\leq 1$, then $\dim_I(X)\leq 1$ if and only if $\dim_G(X)\leq 1$.
\end{Lem}
\dokaz In view of \ref{0dim} the only interesting case is that
of $\dim_JX=1$.
Use fibration $K(J,1)\to K(G,1)\to K(I,1)$
and the fact $X\tau K(J,1)$ to conclude $X\tau K(I,1)$ if and only if
$X\tau K(G,1)$. 
\hfill $\blacksquare$

If $G$ is a group, then $\Ab(G)$ is its {\it abelianization}.

\begin{Lem} \label{GroupsAndAbelianizations}
If $X$ is a metrizable space, then $\dim_{\Ab(G)}(X)\leq \dim_GX$
for any group $G$.
\end{Lem}
\dokaz In view of \ref{0dim} the only interesting case is that
of $G$ non-Abelian and $\dim_GX=1$. 
Since $X\tau K(G,1)$ one gets $X\tau K(H_1(K(G,1)),1)$ by Theorem 3.4
of  \cite{Dy3}. As $H_1(K(G,1))=\Ab(G)$, we are done.
\hfill $\blacksquare$

\section{Nilpotent groups}

If $A,B\subset G$ are subgroups, then the {\it commutator subgroup} $[A,B]$
is a group generated by all commutators $[a,b]$, $a\in A$ and $b\in B$.
The {\it lower central series} $\{\Gamma_n(G)\}$ for a group $G$ is defined as follows:
$\Gamma_1(G)=G$, and $\Gamma_{n+1}(G)=[\Gamma_n(G),G]$.
If a group $G$ is nilpotent, then there exists an integer $c$ such that 
$\Gamma_{c}(G)\ne\{1\}$ but $\Gamma_{c+1}(G)=\{1\}$.
The number $c$ (denoted by $h(G)$) is called the {\it nilpotency class} of the nilpotent group $G$ or
its {\it Hirsch length}. Abelian groups are nilpotent of Hirsch length $1$.
By \cite[Theorem 3.1]{W}, for every $n$ there exists an epimorphism
$$
\otimes^n \Ab(G)\to \Gamma_n(G)/\Gamma_{n+1}(G).
$$
In particular, there is an epimorphism $\otimes^c \Ab(G)\to \Gamma_c(G)$.
It follows from the definition that $\Gamma_{c}(G)$ is in the center of $G$.
Therefore $1\to \Gamma_{c}(G)\to G \to G/\Gamma_{c}(G)\to 1$ is a central extension.
A short calculation shows that the epimorphism $G \to G/\Gamma_{c}(G)$ induces the
trivial homomorphism $\Gamma_c(G)\to \Gamma_c(G/\Gamma_{c}(G))=1$. Therefore 
$G/\Gamma_{c}(G)$ is nilpotent of nilpotency class strictly less than $c$.
This motivates the following definition.

\begin{Def}
A central extension $K\to G\to I$ of groups where $G$ is nilpotent (or equivalently $I$ is nilpotent), 
for which there exists an epimorphism $\otimes^n \Ab(G)\to K$ for some $n$, 
is called a {\it nilpotent central extension}.
\end{Def}

Thus, for every (nonabelian) nilpotent group $G$, there exists a nilpotent central extension
$K\to G\to I$ such that the Hirsch length of $I$ is less than the Hirsch length of $G$.

\begin{Lem}\label{pdivisibility}
Let $1 \to K \to G \stackrel{\pi}{\to} I \to 1$ be a central extension
of nilpotent groups. 
\begin{itemize}
\item[(a)] If $K$ and $I$ are $p$-divisible then $G$ is $p$-divisible.
\item[(b)] If the extension is a nilpotent central extension and $G$ is $p$-divisible, 
then $K$ and $I$ are $p$-divisible.
\end{itemize}
\end{Lem}
\dokaz
Suppose $K$ and $I$ are $p$-divisible.
Let $g\in G$. Then $\pi(g)=i^p$ for some $i\in I$. Let $\bar g\in G$ be such
that $\pi(\bar g)=i$. Then $\pi({\bar g}^p g^{-1})=1$, so ${\bar g}^p g^{-1}=k^p$
for some $k\in K$. Because $k\in K\subset C(G)$, $g=(\bar g k^{-1})^p$, so $G$ is $p$-divisible.

If $G$ is $p$-divisible, then any epimorphic image of $G$ is $p$-divisible. Thus both $I$
and $\Ab(G)$ are $p$-divisible. As there is an epimorphism
$\otimes^n \Ab(G)\to K$, $K$ is also $p$-divisible.
 \hfill
$\blacksquare$

\begin{Lem}\label{ComparisonOfNilpotentClasses}
Suppose $\PPP_i$, $i=1,2$, are two classes of nilpotent groups
such that for any nilpotent central extension
$1\to K\to G\to I\to 1$ where $h(I)<h(G)$ the following conditions hold
\begin{itemize}
\item[(a)] $K$ and $I$ belong to $\PPP_1$ if $G\in\PPP_1$,
\item[(b)] $G\in\PPP_2$ if $K,I\in \PPP_2$.
\end{itemize}
If $A\in\PPP_1\implies A\in\PPP_2$ for all Abelian groups $A$,
then $\PPP_1\subset \PPP_2$.
\end{Lem}
\dokaz
We prove the implication $G\in\PPP_1\implies G\in\PPP_2$ by induction on the Hirsch length of $G$. 
By assumption the implication holds for Abelian group $G$, therefore suppose 
$1\to K\to G\to I\to 1$ is a nilpotent central extension of groups
and $I$ is of lower Hirsch length than $G$. By (a) we have $K,I\in\PPP_1$. By inductive hypothesis
$K,I\in\PPP_2$. The condition (b) yields $G\in\PPP_2$.
 \hfill
$\blacksquare$

\begin{Cor}\label{HomologyOfNilpotent}
Let $A$ be an Abelian group and $n\in\{1, 2\}$. Consider the following statements:
\begin{enumerate}
\item
$\widetilde{H}_*(G;A)=0$, 
 \item
$\widetilde{H}_i(G;A)=0$ for all $i\leq n$.
\end{enumerate}
If $(1)$ is equivalent to $(2)$ for all Abelian groups $G$,
then the two statements are equivalent for all nilpotent groups $G$.
\end{Cor}
\dokaz
Let $\PPP_1^r$ (respectively, $\PPP_2^r$) be the class of all nilpotent
groups $G$ such that $\widetilde{H}_i(G;A)=0$ for all $i\leq n$ (respectively,
$\widetilde{H}_i(G;A)=0$ for all $i$) and $h(G)\leq r$.
Our goal is to prove, by induction on $r$, that $\PPP_1^r=\PPP_2^r$.
It is clearly so for $r=1$. Assume $\PPP_1^m=\PPP_2^m$ for all $m < r$.

Suppose $1\to K\to G\to I\to 1$ is a nilpotent central extension of groups
such that $h(I)< h(G)=r$.
If $G\in \PPP_j^r$, then $H_1(K;A)=0$
which implies $H_i(G;A)\to H_i(I;A)$ is an
epimorphism for $i\leq 2$, so
$H_i(I;A)=0$ for $1\leq i\leq n$. By inductive assumption
$H_i(I;A)=0$ for all $i\ge 1$. If $\widetilde{H}_i(K;A)$ is the first non-trivial reduced homology group of $K$, then Leray-Serre spectral sequence implies
$H_{i+1}(I;A)\ne 0$, a contradiction. Thus both $K$ and $I$
belong to $\PPP_j^r$.
Conversely, if $K,I\in\PPP_j^{r-1}$, then (by inductive assumption) they
have trivial homology with coefficients in $A$ resulting in $G$ having trivial homology with coefficients in $A$ and $G\in\PPP_j^r$.
Applying \ref{ComparisonOfNilpotentClasses} one gets $\PPP_1^r=\PPP_2^r$.
\hfill
$\blacksquare$

\section{Bockstein basis}

If $G$ is a group, then  $\Tor (G)$ is the
subgroup generated by torsion elements of $G$, $\Tor_p(G)$ is the subgroup generated
by all elements of $G$ whose order is a power of $p$, $F_p(G) = G/\Tor_p(G)$, and
$F(G) = G/\Tor(G)$. 

The {\it Bockstein groups} are: rationals $\Q$, cyclic
groups $\Z/p$ of $p$ elements, $p$-adic circles $\Z/p^\infty$, and
$p$-localizations of integers $\Z_{(p)}=\{\frac m n\in\Q\mid n\;\text{is not divisible by}\; p\}$,
where $p$ is a prime number.
Here is a classical definition of the maximal Bockstein basis of an Abelian group
$G$:

\begin{Def}\label{DefOfMaxBocksteinBasis}
$\sigma (G)$ is a subset of all Bockstein groups satisfying the following conditions:
\begin{enumerate}
\item $\Q\in \sigma(G)$ if and only if $F(G)\not= 1$.
\item $\Z/p^\infty\in\sigma(G)$ if and only if $\Tor_p(G)\ne 1$ or
$F(G)$ is not divisible by $p$.
\item $\Z/p\in\sigma(G)$ if and only if $\Tor_p(G)$ is not divisible by $p$ or $F(G)$ is not divisible by $p$.
\item $\Z_{(p)} \in \sigma(G)$ if and only if $F(G)$ not divisible by $p$.
\end{enumerate}
\end{Def}

Notice that our definition \ref{DefOfMaxBocksteinBasis} differs from that in \cite{DR4} 
in the sense that ours is maximal (if $\Z/p\in\sigma(G)$,
then $\Z/p^\infty\in\sigma(G)$ and if $\Z_{(p)} \in \sigma(G)$, then $\Z/p\in\sigma(G)$)
and the one in \cite{DR4} is minimal (if $\Z/p\in\sigma(G)$ or $\Z_{(p)} \in \sigma(G)$, then $\Z/p^\infty\notin\sigma(G)$). From the point of view of First Bockstein Theorem both definitions
are equivalent.

Here is a definition for nilpotent groups which is more convenient in this paper
as it allows using localization of short exact sequences of nilpotent groups.
Recall that a (multiplicative) group is $p$-local iff the map $x\mapsto x^p$ is a bijection.
We call a nilpotent group $\bar p$-local iff it is $q$-local for all prime $q\ne p$.

\begin{Def}\label{BocksteinBasisForNil}
Let $G$ be a nilpotent group, then the Bockstein basis $\sigma(G)$ 
is defined as follows:
\begin{enumerate}
\item $\Q \not\in \sigma(G)$ if and only if $G = \Tor(G)$.
\item $\Z/p^\infty\not\in\sigma(G)$ if and only if $G$ is $\bar p$-local.
\item $\Z/p\not\in\sigma(G)$ if and only if $G$ is divisible by $p$.
\item $\Z_{(p)}\not\in\sigma(G)$ if and only if $F(G)$ is $\bar p$-local.
\end{enumerate}
\end{Def}

\begin{Rem} \label{InclusionsOfSigma}
Note that according to the above definition we have
$$
\Z_{(p)}\in\sigma(G) \Rightarrow \Z/p\in\sigma(G)\Rightarrow \Z/p^\infty\in\sigma(G).
$$
\end{Rem}

\begin{Cor}\label{CharOfQinBB}
For a nilpotent group $G$ the following statements are equivalent:
\begin{enumerate}
\item $\Q \not\in \sigma(G)$,
\item $\widetilde{H}_*(G;\Q)=0$, and
\item $H_1(G;\Q)=0$.
\end{enumerate}
\end{Cor}
\dokaz
For Abelian groups (2) and (3) are equivalent \cite[Theorem 6.4(ii)]{B}.
With Corollary \ref{HomologyOfNilpotent} we extend this equivalence to all 
nilpotent groups.

Let $\PPP_1$ be the class of torsion nilpotent groups
and let $\PPP_2$ be the class of nilpotent groups such that
$H_1(G;\Q)=0$ (i.e., $\Ab(G)$ is torsion). Use Lemma \ref{ComparisonOfNilpotentClasses}
to conclude $\PPP_1=\PPP_2$ and thus (1) and (3) are equivalent.
\hfill
$\blacksquare$

\begin{Cor}\label{CharOfZ/pinBB}
For a nilpotent group $G$, $\Z/p \not\in \sigma(G)$ if and only if $H_1(G;\Z/p)=0$.
\end{Cor}
\dokaz
Let $\PPP_1$ be the class of $p$-divisible nilpotent groups
and let $\PPP_2$ be the class of nilpotent groups such that
$H_1(G;\Z/p)=0$. Use \ref{ComparisonOfNilpotentClasses}
and \ref{pdivisibility}
to conclude $\PPP_1=\PPP_2$.
\hfill
$\blacksquare$

\begin{Cor} \label{ZpInfty}
For a nilpotent group $G$ the following statements are equivalent:
\begin{enumerate}
\item $\Z/p^\infty \not\in \sigma(G)$, \item
$\widetilde{H}_*(G;\Z/p^\infty)=0$, and \item
$H_1(G;\Z/p^\infty)=H_2(G;\Z/p^\infty)=0$.
\end{enumerate}
\end{Cor}
\dokaz If $\Z/p^\infty\notin \sigma(G)$, then $G$
is $\bar p$-local, so all its integral homology groups are $\bar
p$-local and $H_i(G;\Z/p^\infty)=0$ for all $i\ge 1$, which proves
the implication $(1)\implies (2)$.

Notice $(2)\iff (3)$ by \ref{HomologyOfNilpotent}. Indeed,
if $A$ is Abelian and $H_1(A;\Z/p^\infty)=H_2(A;\Z/p^\infty)=0$,
then $A$ must be $\bar p$-local.

Let $\PPP_1$ be the class of $\bar p$-local nilpotent groups
and let $\PPP_2$ be the class of nilpotent groups such that
$H_1(G;\Z/p^\infty)=H_2(G;\Z/p^\infty)=0$. Use \ref{ComparisonOfNilpotentClasses}
to conclude $\PPP_1=\PPP_2$.
\phantom{aaaaaaaaaaaaaaaaaaa}\hfill
$\blacksquare$

\begin{Cor}\label{CharZLocAtpinBB}
For a nilpotent group $G$ the following statements are equivalent:
\begin{enumerate}
\item $\Z_{(p)} \not\in \sigma(G)$, \item
$\widetilde{H}_*(F(G);\Z/p^\infty)=0$, and \item
$H_1(F(G);\Z/p^\infty)=H_2(F(G);\Z/p^\infty)=0$.
\end{enumerate}
\end{Cor}

It is obvious that $\Z/p^\infty \not\in \sigma(G)$ if and only if $G\to G_{(\bar p)}$ is an isomorphism.
The following lemma characterizes $\Z_{(p)} \not\in \sigma(G)$ via localizations.
\begin{Lem}\label{CharZLocAtpinBBbyLocalizations}
For a nilpotent group $G$ the following statements are equivalent:
\begin{enumerate}
\item $\Z_{(p)} \not\in \sigma(G)$, 
\item $G\to G_{(\bar p)}$ is an epimorphism.
\end{enumerate}
\end{Lem}
\dokaz If $G\to G_{(\bar p)}$ is an epimorphism, then its kernel has trivial
$\bar p$-localization (by exactness of the localization
functor) and must be a torsion group. Therefore $F(G)=F(G_{(\bar p)})$ is $\bar p$-local.
If $F(G)$ is $\bar p$-local, then apply exactness of the localization functor
to the short exact sequence $1\to \Tor(G)\to G\to F(G)\to 1$
and derive $G\to G_{(\bar p)}$ is an epimorphism.
\hfill
$\blacksquare$

\begin{Def}
The {\it torsion-divisible Bockstein basis} $\sigma_{TD}(G)$ of $G$ consists of all
$\Z/p^\infty$ belonging to $\sigma(G)$. We set $\sigma_{NTD}(G) = \sigma(G) \backslash \sigma_{TD}(G)$.
\end{Def}

\begin{Lem}\label{EpimorphicImage}
If $G\to I$ is an epimorphism of nilpotent groups, then
$\sigma_{NTD}(I)\subset\sigma_{NTD}(G)$.
\end{Lem}
\dokaz
Suppose $\Q\not\in\sigma(G)$, then $G$ is a torsion group.
So $I$ is a torsion group, hence $\Q\not\in\sigma(I)$.

Let $\Z/p\not\in\sigma(G)$, then $G$ is $p$-divisible and then
also $I$ is $p$-divisible.

Let $\Z_{(p)}\not\in\sigma(G)$, then $G\to G_{(\bar p)}$ is an epimorphism.
Because $\bar p$-localization is an exact functor, the map
$G_{(\bar p)}\to I_{(\bar p)}$ is an epimorphism and hence
$I\to I_{(\bar p)}$ is an epimorphism.
 \hfill
$\blacksquare$

\begin{Lem} \label{SigmaGsubsetSigmaKandSigmaI}
Let $1 \to K \to G \to I \to 1$ be a central extension of nilpotent groups,
then $\sigma(G)\subset\sigma(K)\cup\sigma(I)$.
\end{Lem}
\dokaz
Let $\Q\not\in \sigma(K)\cup\sigma(I)$, then $K$ and $I$ are torsion groups.
Hence $G$ is torsion, so $\Q\not\in\sigma(G)$.

Let $\Z/p\not\in \sigma(K)\cup\sigma(I)$, then $K$ and $I$ are $p$-divisible.
By Lemma \ref{pdivisibility}, also $G$ is $p$-divisible.

Let $\Z/p^\infty\not\in \sigma(K)\cup\sigma(I)$, then
$K\to K_{(\bar p)}$ and $I\to I_{(\bar p)}$ are isomorphisms.
Using Five Lemma and the fact that $\bar p$-localization is
an exact functor, we conclude that also $G\to G_{(\bar p)}$ is an isomorphism.

Let $\Z_{(p)}\not\in \sigma(K)\cup\sigma(I)$, then
$K\to K_{(\bar p)}$ and $I\to I_{(\bar p)}$ are epimorphisms.
By Three Lemma \cite[Lemma 2.8]{H}, also $G\to G_{(\bar p)}$ is an epimorphism.
 \hfill
$\blacksquare$

\begin{Lem} \label{SigmaAbGsubsetSigmaAbKandSigmaAbI}
Let $1 \to K \to G \to I \to 1$ be a nilpotent central extension of groups.
If $\Z_{(p)}\not\in\sigma(\Ab(G))$ for some prime $p$,
then
$\Z_{(p)}\not\in\sigma(K)$ and $\Z_{(p)}\not\in\sigma(\Ab(I))$.
\end{Lem}
\dokaz
Assume $\Z_{(p)}\not\in\sigma(\Ab(G))$. That means the map $F(\Ab(G))\to F(\Ab(G))_{(\bar p)}$ is an isomorphism.
Because $\Ab(G)\to \Ab(I)$ is an epimorphism and $F$ is a right exact functor,
the map $F(\Ab(G))\to F(\Ab(I))$ is an epimorphism.
Hence the map $F(\Ab(I))\to F(\Ab(I))_{(\bar p)}$ is an epimorphism. Its kernel is
a $p$-torsion group, so the kernel is trivial and the map
$F(\Ab(I))\to F(\Ab(I))_{(\bar p)}$ is an isomorphism.
That means $\Z_{(p)}\not\in\sigma(\Ab(I))$.

There exists an epimorphism $\otimes^n \Ab(G)\to K$ for some integer $n$.
Because $F(\Ab(\otimes^n G))\to F(\Ab(\otimes^n G))_{(\bar p)}$ is an isomorphism,
in the same way as in the previous paragraph we can
prove that $F(K)\to F(K)_{(\bar p)}$ is an isomorphism.
Hence $\Z_{(p)}\not\in\sigma(K)$.
 \hfill
$\blacksquare$

\begin{Lem} \label{SigmaAbGsubsetSigmaG}
Let $G$ be a nilpotent group, then $\sigma_{NTD}(G)=\sigma_{NTD}(\Ab(G))$
and $\sigma(\Ab(G))\subset\sigma(G)$.
\end{Lem}
\dokaz
The inclusion $\sigma_{NTD}(\Ab(G))\subset \sigma_{NTD}(G)$ follows from Lemma \ref{EpimorphicImage}.

Let us prove $\sigma_{NTD}(G)\subset\sigma_{NTD}(\Ab(G))$.
Suppose $\Q\not\in\sigma_{NTD}(\Ab(G))$. Because $G$ is a torsion group if and only if $\Ab(G)$
is a torsion group, $\Q\not\in\sigma_{NTD}(G)$.

Suppose $\Z/p\not\in\sigma_{NTD}(\Ab(G))$. Because $G$ is $p$-divisible if and only if $\Ab(G)$
is $p$-divisible, $\Z/p\not\in\sigma_{NTD}(G)$.

Consider the class $\PPP$ of all nilpotent groups such that $\Z_{(p)}\not\in\sigma_{NTD}(\Ab(G))$ implies
$\Z_{(p)}\not\in\sigma_{NTD}(G)$. $\PPP$ clearly
contains all Abelian groups. To show $\PPP$ equals the class $\NNN$ of all nilpotent groups
it suffices to show (see \ref{ComparisonOfNilpotentClasses}) that for any nilpotent central extension $A\to G\to G'$ such that Hirsch length of $G'$ is less
than $h(G)$,  $A,G'\in\PPP$ implies $G\in \PPP$.
Assume $\Z_{(p)}\not\in\sigma_{NTD}(\Ab(G))$. 
By \ref{SigmaAbGsubsetSigmaAbKandSigmaAbI} we conclude $\Z_{(p)}\not\in\sigma_{NTD}(G')$ as $G'\in\PPP$ and $\Z_{(p)}\not\in\sigma_{NTD}(A)$.
By \ref{SigmaGsubsetSigmaKandSigmaI} $\Z_{(p)}\not\in\sigma_{NTD}(G)$.

Let us prove now that $\sigma(\Ab(G))\subset\sigma(G)$. By Lemma \ref{EpimorphicImage},
$\sigma_{NTD}(\Ab(G))\subset\sigma_{NTD}(G)$.
Suppose $\Z/p^\infty\not\in \sigma(G)$. Then $G$ is uniquely $p$-divisible,
so $H_*(G;\Z)$ is uniquely $p$-divisible. In particular $\Ab(G)=H_1(G;\Z)$
is uniquely $p$-divisible, hence $\Z/p^\infty\not\in \sigma(\Ab(G))$.
 \hfill
$\blacksquare$

\begin{Thm} \label{SigmaG=SigmaKandSigmaI}
If $1 \to K \to G \to I \to 1$ is a nilpotent central extension, then $\sigma(G) =
\sigma(K) \cup \sigma(I)$.
\end{Thm}
\dokaz
By Lemma \ref{SigmaGsubsetSigmaKandSigmaI}, $\sigma(G) \subset\sigma(K) \cup \sigma(I)$.

Let us prove that $\sigma(K) \cup \sigma(I)\subset \sigma(G)$.
Suppose $\Q\not\in\sigma(G)$, then $G$ is a torsion group. Therefore
$K$ and $I$ are also torsion groups, so $\Q\not\in \sigma(K) \cup \sigma(I)$.

Suppose $\Z/p\not\in\sigma(G)$, then $G$ is $p$-divisible. By Lemma \ref{pdivisibility},
$K$ and $I$ are $p$-divisible, so $\Z/p\not\in \sigma(K) \cup \sigma(I)$.

Suppose $\Z/p^\infty\not\in\sigma(G)$, then $G\to G_{(\bar p)}$ is an isomorphism.
Because $\bar p$-localization is an exact functor, the map
$K\to K_{(\bar p)}$ is a monomorphism and the map $\Ab(G)\to (\Ab(G))_{(\bar p)}$ is
an epimorphism. Because there exists an epimorphism $\otimes^n \Ab G\to K$,
the map $K\to K_{(\bar p)}$ is also an epimorphism, hence it is an isomorphism.
By Five Lemma, also the map $I\to I_{(\bar p)}$ is an isomorphism,
so $\Z/p^\infty\not\in \sigma(K) \cup \sigma(I)$.

If $\Z_{(p)}\not\in\sigma(G)$, then $\Z_{(p)}\not\in\sigma(K)$
and $\Z_{(p)}\not\in\sigma(I)$ by \ref{SigmaAbGsubsetSigmaAbKandSigmaAbI}
and \ref{SigmaAbGsubsetSigmaG}.
 \hfill
$\blacksquare$

\section{Bockstein spaces}

\begin{Def}
A metrizable space $X$ is called a {\it Bockstein space} if 
$\dim_G X=\sup\{\dim_H X\mid H\in\sigma(G)\}$ for all Abelian groups $G$.
\end{Def}

\begin{Rem}
In the above definition observe $\dim_G X$ is an element of $\N\cup\{0,\infty\}$ and not only
in $\{0,1,\infty\}$ as in the case of non-Abelian groups $G$.
\end{Rem}

Dranishnikov-Repov\v s-Shchepin \cite{DrRS}
showed the existence of a separable metric space $X$ of dimension $2$ such that
$\dim_{\Z_{(p)}}X=1$ for all primes $p$. Thus, $X$ is not a Bockstein space
as $\dim_\Z X=2 > 1=\sup\{\dim_H X\mid H\in\sigma(\Z)\}$.

\begin{Problem}\label{AreANRSBockstein}
Is every metric ANR a Bockstein space?
\end{Problem}

\begin{Prop}\label{FSigmaIsBockstein}
Suppose $X=\bigcup\limits_{n=1}^\infty X_n$ is metrizable
and each $X_n$ is closed in $X$. If all $X_n$ are Bockstein spaces,
then so is $X$.
\end{Prop}
\dokaz Suppose $G$ is an Abelian group and $H\in\sigma(G)$.
If $\dim_G(X)\leq m$, then $\dim_G(X_n)\leq m$ for all $n$
and $\dim_H(X_n)\leq m$ for all $n$ resulting in $\dim_H(X)\leq m$.
\par If $\dim_H(X) < m$ for all $H\in \sigma(G)$, then 
$\dim_H(X_n) < m$ for all $n$ and $\dim_G(X_n) < m$ for all $n$
resulting in $\dim_G(X) < m$.
\hfill
$\blacksquare$

For a subset $l\subset \mathbb P$ of prime integers let
$\Z_l=\{\frac m n\in\Q\mid n\;\text{is not divisible by any}\; p\in l  \}$ and
let $\hat\Z_l$ be the group of $l$-adic integers.

\begin{Thm}\label{CharOfBocksteinSpaces}
A metrizable space $X$ is a Bockstein space if and only if $\dim_{\Z_l}X=\dim_{\hat\Z_l} X$
for all subsets $l\subset \mathbb P$ of the set of prime numbers.
\end{Thm}
\dokaz
Since $\sigma(\hat\Z_l)=\sigma(\Z_l)$ for all $l\subset \mathbb P$,
$\dim_{\Z_l}X=\dim_{\hat\Z_l} X$ holds for any Bockstein space $X$.
\par Assume $\dim_{\Z_l}X=\dim_{\hat\Z_l} X$
for all subsets $l\subset \mathbb P$ of the set of prime numbers.
Suppose $G$ is a torsion-free Abelian group $G$.
If $\Z_{(p)}\in \sigma(G)$, then Theorem B(d) of \cite{Dy4}
says $\dim_{\hat\Z_{(p)}}X\leq \dim_GX$. Therefore $\dim_G X\ge\sup\{\dim_H X\mid H\in\sigma(G)\}$.
Suppose $\sup\{\dim_H X\mid H\in\sigma(G)\}=n$ and consider $l=\{p\mid p\cdot G\ne G\}$.
Theorem B(f) of \cite{Dy4} says $\dim_G(X)\leq \dim_{\Z_l}X$.
Since $\sigma(G)=\sigma(\Z_l)$, $\dim_GX\leq \dim_{\Z_l}X=
\sup\{\dim_H X\mid H\in\sigma(\Z_l)\}=\sup\{\dim_H X\mid H\in\sigma(G)\}=n$.
That proves $\dim_G X=\sup\{\dim_H X\mid H\in\sigma(G)\}$
for all torsion-free Abelian groups. The same equality holds for all torsion Abelian groups by Theorem B(a) of \cite{Dy4}.
In the case of arbitrary Abelian groups $G$, as $\sigma(G)=\sigma(F(G))\cup \sigma(\Tor(G))$
and $\dim_G=\max(\dim_{F(G)}X,\dim_{\Tor(G)}X)$ (see Theorem B(b) of \cite{Dy4})
one gets $\dim_G X=\sup\{\dim_H X\mid H\in\sigma(G)\}$ as well.
 \hfill
$\blacksquare$

\begin{Rem}
Notice that it is not sufficient to assume $\dim_{\Z_{(p)}}X=\dim_{\hat\Z_{(p)}} X$
for all primes $p$ in
\ref{CharOfBocksteinSpaces}. Indeed, the space $X$ in \cite{DrRS}
has that property as
$1=\dim_{\Z_{(p)}}X\ge \dim_{\hat\Z_{(p)}} X\ge 1$ for all primes $p$. 
\end{Rem}

\begin{Thm}\label{BocksteinForNilpotentThm}
Let $X$ be a Bockstein space.
If $G$ is nilpotent,
then $\dim_G(X) \leq 1$ if and only if $\sup\{\dim_H(X) | H\in\sigma(G)\}\leq 1$.
\end{Thm}
\dokaz
 \ref{GroupsAndAbelianizations}
Let $\PPP_1$ be the class of all nilpotent groups and let $\PPP_2$ be the class of nilpotent
groups $G$ such that $\dim_G(X) \leq 1$ if and only if $\sup\{\dim_H(X) | H\in\sigma(G)\}\leq 1$.
Since $\PPP_2$ contains all Abelian groups, in view of \ref{ComparisonOfNilpotentClasses}
it suffices to show that for any nilpotent central extension $K\to G\to I$ the conditions $K,I\in \PPP_2$ imply $G\in \PPP_2$. It is so if $G$ is Abelian, so assume $G$
is not Abelian. Moreover, as $\sigma(\Ab(G))\subset \sigma(G)$ by
\ref{SigmaAbGsubsetSigmaG} and $\dim_GX\leq 1$ implies  $\dim_{\Ab(G)}X\leq 1$
(see  \ref{GroupsAndAbelianizations}),
either $\dim_G(X) \leq 1$ or $\sup\{\dim_H(X) | H\in\sigma(G)\}\leq 1$
implies $\dim_{\Ab(G)}X\leq 1$, so we may as well assume $\dim_{\Ab(G)}X\leq 1$.

In view of Lemma \ref{0dim} and the fact $\Ab(G)=1$ implies $G=1$, the equivalence
of conditions $\dim_G(X) \leq 1$ and $\sup\{\dim_H(X) | H\in\sigma(G)\}\leq 1$
may fail only if $\dim_{\Ab(G)}X=1$, so assume $\dim_{\Ab(G)}X=1$.
\par
Suppose $\dim_H(X)=n > 1$ for some $H\in\sigma(K)$. If $H\ne\Z/p^\infty$, then
$H\in\sigma_{NTD}(G)=\sigma_{NTD}(\Ab(G))$ (Theorem \ref{SigmaG=SigmaKandSigmaI} and Lemma \ref{SigmaAbGsubsetSigmaG}),
a contradiction. So $H=\Z/p^\infty$ for some prime $p$ and $\Z/p^\infty\not\in\sigma(\Ab(G))$.
Hence
$\Ab(G)$ is $p$-divisible and this is equivalent to $G$ being 
$p$-divisible \cite[Lemma 5.1]{CDMV}.
Because $\Z/p^\infty\not\in\sigma(\Ab(G))$, by Lemma \ref{ZpInfty}
$H_1(\Ab(G); \Z/p^\infty)=0$, so also $H_1(G; \Z/p^\infty)=0$.
By Lemma \ref{ZpInfty}, $H_2(G; \Z/p^\infty)\ne 0$ as $\Z/p^\infty\in \sigma(G)$.
This implies $F(H_2(G; \Z))$ is not $p$-divisible, so $\Z_{(p)}\in \sigma(H_2(G; \Z))$.
If $G$ is not a torsion group, then
 $\Ab(G)$ is not a torsion group, hence by definition $\Q\in\sigma(\Ab(G))$.
Therefore $\dim_{\Q}(X)\leq 1$. Using that fact and
Bockstein Inequalities (BI5, BI6 \cite{K}),
we get $\dim_{\Z/p^\infty}(X)=\dim_{\Z_{(p)}}(X)-1$.
Because $\Z_{(p)}\in\sigma(H_2(G;\Z))$, the dimension
$\dim_{\Z_{(p)}}(X)\le\dim_{H_2(G;\Z)}(X)\le 2$ as $X$ is a Bockstein space and then
$\dim_{\Z/p^\infty}(X)=\dim_{\Z_{(p)}}(X)-1\le 1$, a contradiction.

Thus $G$ is a torsion group and is a product
of $q$-groups $G=\Pi_{q\in\PP} G_q$. Hence $\Ab(G)=\Pi_{q\in\PP} \Ab(G_q)$.
Because $G$ is not $\bar p$-local, $G_p\ne 1$, but $\Ab(G)$ is uniquely $p$-divisible,
so $\Ab(G_p)=1$. Therefore $G_p$ is a perfect nilpotent group, 
but such group is trivial, a contradiction.

Thus $\dim_HX\leq 1$ for all $H\in\sigma(K)$ and
$\dim_KX=\sup\{\dim_H(X) | H\in\sigma(K)\}\leq 1$ as $K$
is Abelian and $X$ is a Bockstein space.

Suppose $\sup\{\dim_H(X) | H\in\sigma(G)\}\leq 1$.
 By Theorem \ref{SigmaG=SigmaKandSigmaI},
$\sigma(G) =\sigma(K) \cup \sigma(I)$. Therefore $\sup\{\dim_H(X) | H\in\sigma(I)\}\leq 1$
and $\dim_IX\leq 1$ as $I\in\PPP_2$. Consequently,
 $\dim_GX\leq1$ by \ref{BocksteinInequalities}.
 
 Suppose $\dim_GX\leq 1$. By \ref{BocksteinInequalities}
 one gets $\dim_IX\leq 1$  and $\sup\{\dim_H(X) | H\in\sigma(I)\}\leq 1$
 as $I\in\PPP_2$.
 By Theorem \ref{SigmaG=SigmaKandSigmaI},
$\sigma(G) =\sigma(K) \cup \sigma(I)$, hence
$$
\sup\{\dim_H(X) | H\in\sigma(G)\}=\sup\{\dim_H(X) | H\in\sigma(K) \cup \sigma(I)\}\leq
1.$$ \hfill
$\blacksquare$

\begin{Cor}\label{HurewiczSerre}
Let $L$ be a connected nilpotent CW complex.
If $X$ is a Bockstein space and $\dim_{H_n(L)}X\leq n$ for all $n\ge 1$,
then $\dim_{\pi_n(L)}X\leq n$ for all $n\ge 1$.
\end{Cor}
\dokaz It is shown in \cite{CDMV} that $\dim_{\pi_n(L)}X\leq n$ for all $n\ge 2$,
so it suffices to prove $\dim_{\pi_1(L)}X\leq 1$.
If that inequality is false, then there is $H\in\sigma(\pi_1(L))$
such that $\dim_H(X) > 1$ (see \ref{BocksteinForNilpotentThm}).
In view of \ref{SigmaAbGsubsetSigmaG}, as $\dim_{H_1(L)}(X)\leq 1$,
$H=\Z/p^\infty$ for some prime $p$.
Also, $H_1(L)$ is not a torsion group, so $\dim_{\Q}(X)\leq 1$.
Using Bockstein Inequalities one gets $\dim_{\Z_{(p)}}(X)\ge 3$.
Therefore $i$-th homology groups of both $L$ and $\tilde L$ with coefficients
in $\Z/p^\infty$ vanish for $i\leq 2$. From the fibration
$\tilde L\to L\to K(\pi_1(L),1)$ one gets $i$-th homology groups of $K(\pi_1(L),1)$ with coefficients
in $\Z/p^\infty$ vanish for $i\leq 2$. However, in view of \ref{ZpInfty},
that means $\Z/p^\infty\notin\sigma(\pi_1(L))$, a contradiction.
\hfill
$\blacksquare$

\begin{Cor}\label{HurewiczSerreForFiniteDim}
Let $L$ be a connected nilpotent CW complex
and let $X$ be a Bockstein space such that $\dim_{H_n(L)}X\leq n$ for all $n\ge 1$.
If $X$ is finite dimensional or $X\in ANR$, then $X\tau L$.
\end{Cor}
\dokaz By \ref{HurewiczSerre} one gets $\dim_{\pi_n(L)}X\leq n$ for all $n\ge 1$
and by Theorem G in \cite{Dy4} one has $X\tau L$.
\hfill
$\blacksquare$

\end{document}